\documentclass{ifacconf}
\usepackage{amsmath}
\usepackage{amssymb}
\usepackage{amsfonts}
\usepackage{graphicx}      % include this line if your document contains figures
\usepackage{natbib}        % required for bibliography
\usepackage{url}

\usepackage{booktabs}
\usepackage{color}
\allowdisplaybreaks
\newcommand{\cat}{{\mathrm{cat}}}
\newcommand{\tp}{t_{\mathrm{p}}}
\newcommand{\tc}{t_{\mathrm{c}}}
\newcommand{\ts}{t_{\mathrm{s}}}
\newcommand{\tend}{t_{\mathrm{end}}}
\newcommand{\vlin}{v_{\mathrm{lin}}}
\newcommand{\bal}{{\mathrm{bal}}}
\newcommand{\lin}{{\mathrm{lin}}}
\newcommand{\init}{{\mathrm{init}}}
\newcommand{\mix}{{\mathrm{mix}}}
\newcommand{\opt}{{\mathrm{opt}}}
\newcommand*\diff{\mathop{}\!\mathrm{d}}
\newtheorem{assumption}{Assumption}
\newtheorem{theorem}{Theorem}
\newtheorem{remark}{Remark}

\begin{document}
\begin{frontmatter}

\title{Dynamic modeling of enzyme controlled metabolic networks using a receding time horizon} 

\thanks[footnoteinfo]{H.L. and A.-M.R are funded by ERANET for Systems Biology ERASysApp, project ROBUSTYEAST, BMBF grant
    IDs 031L0017A and 031L0017B.}

\author[First]{Henning Lindhorst} 
\author[Second]{Alexandra-M. Reimers} 
\author[Third]{Steffen Waldherr}

\address[First]{Institute for Automation Engineering, Otto-von-Guericke-Universit\"{a}t Magdeburg; (e-mail: henning.lindhorst@ovgu.de).}
\address[Second]{Department of Mathematics and Computer Science, Freie Universit\"{a}t Berlin; (email: alexandra.reimers@fu-berlin.de)}
\address[Third]{KU Leuven, Department of Chemical Engineering; (e-mail: steffen.waldherr@kuleuven.be)}

\begin{abstract}                % Abstract of not more than 250 words.
  Microorganisms have developed complex regulatory features controlling their reaction and internal adaptation to changing environments.
  When modeling these organisms we usually do not have full understanding of the regulation and rely on substituting it with an
  optimization problem using a biologically reasonable objective function. 
  The resulting constraint-based methods like the Flux Balance Analysis (FBA)
  and Resource Balance Analysis (RBA) 
  have proven to be powerful tools to predict growth rates, by-products, and pathway usage for fixed environments.
  In this work, we focus on the dynamic enzyme-cost Flux Balance Analysis (deFBA), which models the environment, biomass products, and their composition dynamically and contains reaction rate constraints based on enzyme capacity.
  We extend the original deFBA formalism to include storage molecules and biomass-related maintenance costs.
  Furthermore, we present a novel usage of the receding prediction horizon as used in Model Predictive Control (MPC) in the deFBA framework, which we call the short-term deFBA (sdeFBA).
  This way we eliminate some mathematical artifacts arising from the formulation as an optimization problem and gain access to new applications in MPC schemes.
  A major contribution of this paper is also a systematic approach for choosing the prediction horizon and identifying conditions to ensure solutions grow exponentially.
  We showcase the effects of using the sdeFBA with different horizons through a numerical example.
\end{abstract}


\end{frontmatter}
%===============================================================================

\section{Introduction}
Microorganisms encounter a vast array of environmental conditions and have developed complex regulatory mechanisms to cope with them.
While a lot of research is done to investigate this, most regulatory features are still unknown.
An effective alternative approach is the substitution of the regulation with an optimization problem as originally done with the Flux Balance Analysis (FBA) in \citep{varma1994metabolic}. 
This method models the organism as a metabolic network in steady-state and maximizes a single biomass flux.
This approach led to a family of methods focusing on different aspects.

Initial steps towards dynamic models with the ability to react to changing environments were made with the dynamic FBA \citep{mahadevan2002dynamic}.
But this method still lacks a connection between reaction rates and the enzyme levels necessary to realize them.
%Because they assume a fixed biomass composition, these methods can not model the adaptation process under a changing environment as these result in changes in expression levels of enzymes and a 
%The reaction to nutritial changes is mainly done by changing the expression levels of enzymes and thus up- or downregulating different metabolic pathways.
The first optimization method to take this into account is the Resource Balance Analysis (RBA) \citep{goelzer2011cell}.
In this method the growth rate of the cell is optimized to a fixed medium composition while enzymatic flux constraints limit uptake and metabolic reaction rates.
The combination of these enzymatic constraints and a dynamic approach resulted in the dynamic enzyme-cost Flux Balance Analysis (deFBA) presented in \citep{waldherr2015}.
The deFBA predicts all reaction rates and enzymatic levels for given nutrient dynamics on a chosen time frame.
An application of the deFBA to a genome scale model can be found in \citep{reimers2017}.
%One of the benefits of all these methods is the good availability of metabolic models.
%These can then be systematically extended to be used in the deFBA \citep{reimers2017generating}.

During a recent study \citep{waldherr2017} we learned that the fixed end-time in the deFBA can lead to artificial solutions usually not observed in the modeled organisms.
Furthermore, we plan to use deFBA inside a model predictive controller to maximize certain biomass products by manipulation of the medium composition.
Thus, we present in this work the \emph{short-term deFBA} (sdeFBA), which combines the deFBA with the idea of a receding prediction horizon.
This also allows us to solve problems with large end-times piece-wise and in some cases reduces the computational cost for the simulation.

%The work is structured as following.
%First, we present the basics of the deFBA in Section \ref{sec:deFBA} and discuss important growth modes.
%Afterwards, we present the sdeFBA in Section \ref{sec:receding_prediction_horizon} and present in sections \ref{sec:choosing_the_prediction_horizon} and \ref{sec:choosing_the_control_horizon} how to choose a set of prediction and iteration time minimizing the computational cost whilst ensuring qualitative desired results.
%We showcase the effects of using the sdeFBA through a small numerical example in Section \ref{sec:example_1}.

\section{Dynamic enzyme-cost Flux Balance Analysis}\label{sec:deFBA}
\subsection{Constructing the optimization problem}
In this section we present the basics of the deFBA and showcase the extensions of our current formulation in comparison to the original one \citep{waldherr2015}.
At the heart of deFBA models lies a metabolic reaction network consisting of $n$ biochemical species and $m$ reactions converting the species into each other.
We further classify the species depending on their physical location and their biological function as either
\begin{itemize}
    \item \emph{external species} $Y \in \mathbb{R}^{n_y}_{\geq 0}$ outside of the cell (carbon sources, oxygen, etc.),
    \item \emph{metabolic species} $X \in \mathbb{R}^{n_x}_{\geq 0}$ which are intermediates and intracellular products of the metabolism (amino acids, ATP, etc.),
    \item \emph{storage species} $C \in \mathbb{R}^{n_c}_{\geq 0}$ which are allowed to accumulate in the model (glycogen, starch, etc.),
    \item \emph{macromolecules} $P \in \mathbb{R}^{n_p}_{\geq 0}$ representing biomass components (enzymes, cell walls, DNA, etc.),
\end{itemize} 
with $n=n_y + n_x + n_p + n_c$.
We measure all species in molar amounts, e.g., $[X]=$ mol. %and the time variable $t \geq 0,~[t]=$~h.
%Please note that we include all macromolecules and external species as dynamic variables.

The macromolecules $P$ represent the complete reproductive machinery of the organism and can be further divided into a catalytic part, enabling reactions via enzymes and taking care of reproduction via the ribosome, and a non-catalytic part, like cell walls, DNA, etc.
To keep the notation simple we address both kinds with $P$. 
Most organisms use some of the available nutrients to create an energy storage, which can be used to survive phases of starvation, e.g. production of starch during day for consumption at night.
% Plants for example create starch during the day for consumption during the night.
The storage species $C$ can either be some macromolecules or simply metabolites allowed to accumulate.

The deFBA assumes the network maximizes biomass accumulation over time.
Thus, we assign the accumulating species $C, P$ their molecular weights $w_i, ~[w_i]=$ g/mol and define the \emph{total biomass} $B$ as
\begin{align}
    B (t) = w_C^T C(t) + w_P^T P(t),
\end{align}
depending on the time $t$, $[t] =$ h.
As recent studies have shown \citep{waldherr2017} the inclusion of non-catalytic biomass in the objective may lead to unexpected results if these species are very "cheap" to produce in comparison to their weights $w_C$.
Thus, we additionally define the \emph{objective biomass} $B_o$ via the \emph{objective weights} $b_i$, which in most cases coincide with the molecular weights, but can be set to zero if necessary
\begin{align}\label{eq:objective_biomass}
B_o(t) =  b_C^T C(t) + b_P^T P(t).
\end{align}

The reactions $R$ between the species are subdivided into the following types:
\begin{itemize}
    \item \emph{exchange reactions} $v_Y \in \mathbb{R}^{m_y}$ exchanging matter with the outside,
    \item \emph{metabolic reactions} $v_X \in \mathbb{R}^{m_x}$ transforming metabolites into one another,
    \item \emph{storage reactions} $v_C \in \mathbb{R}^{m_c}$ converting metabolites in storage and vice versa,
    \item \emph{biomass reactions} $v_P \in \mathbb{R}^{m_p}$ producing macromolecules, 
\end{itemize}
with $m = m_y + m_x + m_c + m_p$.
We write shortly $v = (v_Y^T, v_X^T, v_C^T, v_P^T)^T$, $[v] = $ mol/h.
The dynamics of the species are then given by the \emph{stoichiometric matrix} $S \in \mathbb{R}^{n,m}$
\begin{align}\label{eq:full_dynamics}
\begin{split}
\frac{\mathrm{d}}{\mathrm{d}t}\begin{pmatrix}
Y(t) \\ X(t) \\ C(t) \\ P(t)
\end{pmatrix} & = \begin{pmatrix}
                                S_{Y,Y} & 0 & 0 & 0 \\
                                S_{X,Y} & S_{X,X} & S_{X,C} & S_{X,P} \\
                                0 & 0 & S_{C,C} & 0 \\
                                0 & 0 & 0 & S_{P,P} \end{pmatrix} \begin{pmatrix} v_Y(t) \\ v_X(t) \\ v_C(t) \\ v_P(t)\end{pmatrix}\\ 
& = \begin{pmatrix} S_Y \\ S_X \\ S_C \\ S_P \end{pmatrix}\begin{pmatrix} v_Y(t) \\ v_X(t) \\ v_C(t) \\ v_P(t) \end{pmatrix} = S v(t),
\end{split}
\end{align}
with the submatrices $S_{I,J} \in \mathbb{R}^{n_I,m_J}, ~I,J \in \{ Y,X,C,P \}$. % and the time variable $t \geq 0$.%, $[t] = $h.
Following \citep{waldherr2015}, the metabolism is modelled to operate in quasi steady-state.
This translates to the constraint
\begin{align}\label{eq:qss}
\begin{split}
\frac{\mathrm{d}}{\mathrm{d}t}X(t) & = S_X v(t) = 0, ~ \forall t \geq 0.  \\
%\Leftrightarrow ~ v(t) &\in \mathrm{Ker}(S_X), ~ \forall t \geq 0.
\end{split}
\end{align}
%The constraint \eqref{eq:qss} lowers the problem size considerably and aligns with the assumption of a growth optimized network.

The enzymatic biomass catalyzes the reactions in the network and the maximal rates are determined by the reaction-specific \emph{catalytic constants} (or turnover numbers) $k_{\cat,\pm j}$, $j \in \{1, \ldots, m\}$, $[k_{\cat,\pm j}]=\mathrm{h}^{-1}$ and the amount of the respective enzyme ${P_i}$.
We differentiate between the \emph{forward value} $k_{\cat, +j}$ and the \emph{backward value} $k_{\cat, -j}$.

The bounds for the reactions rates are given by
\begin{align}
-v_j \leq k_{\cat,- j} {P_i},~ v_j \leq k_{\cat,+ j} {P_i}.
\end{align}
Furthermore, some enzymes are capable of catalyzing multiple reactions, which we describe with the sets
\begin{align}
 \mathrm{cat}({P_i}) = \{ v_j~|~ {P_i} \text{ catalyzes } v_j  \}.
\end{align}
The corresponding constraint with respect to reversibility of the reactions then reads
\begin{align}\label{eq:ecc_base_constraint}
 \sum_{v_j \in \mathrm{cat}({P_i})} \left| \frac{v_j(t)}{k_{\mathrm{cat},\pm j}} \right| \leq {P_i}(t),~ \forall t \geq 0.
\end{align}
%We transform this into a linear constraint by building all possible sign combinations for reversible reactions.
%\textcolor{red}{This could be left out to bring the page number down.}
%As example, consider $P_1$ catalyzing $v_1$ reversibly and $v_2$ irreversibly.
%Then we can write \eqref{eq:ecc_base_constraint} as
%\begin{align}
% H_{c,1} v &= \begin{pmatrix}
%              k_{\mathrm{cat},+1}^{-1} &  k_{\mathrm{cat},+2}^{-1} & 0 & \ldots & 0 \\
%             -k_{\mathrm{cat},-1}^{-1} &  k_{\mathrm{cat},+2}^{-1} & 0 & \ldots & 0
%             \end{pmatrix} v \leq H_{e,1} P,
%\end{align}
%with the filter matrix $H_{e,1}$
%\begin{align}
% H_{e,1} = \begin{pmatrix} 1 & 0 & \ldots & 0 \\ 1 & 0 & \ldots & 0 \end{pmatrix}.
%\end{align}
%The concatenations of the matrices $H_{c,1}$  to $H_{c,n_p}$ and $H_{e,1}$ to $H_{E,n_p}$, are written as $H_c$ and $H_E$, respectively.
We call the matrix form the \emph{enzyme capacity constraint}
\begin{align}\label{eq:ecc}
 H_c v(t) \leq H_{e} P(t), ~ \forall t \geq 0,
\end{align}
with the filter matrix $H_{e}$.
For more detail on the construction of these matrices see \citep{waldherr2015}.
The constraint \eqref{eq:ecc} is the central constraint in deFBA as it limits growth.
%rate $\mu(t)$
%\begin{align}\label{eq:growth_rate}
%\mu(t) = \frac{1}{B(t)}\frac{dB(t)}{dt}
%\end{align}
In regular FBA the growth rate is constrained by biomass independent constraints
\begin{align}\label{eq:box_constraints}
v_{\mathrm{min}} \leq v(t) \leq v_{\max}
\end{align}
derived from measured reaction rates.
Because all reactions can reach arbitrarily large rates given enough enzyme is present (cf.~\eqref{eq:ecc}), we make the following assumption.
\begin{assumption}\label{as:reversibility}
    The biomass independent constraints \eqref{eq:box_constraints} are only used to define the reversibility of the reactions with $v_{\min}, v_{\max} \in \{\pm \infty, 0\}^m$.
\end{assumption}
Any organism needs structural macromolecules to keep working, e.g., the cell wall separating it from the outside.
We express this necessity by enforcing certain fractions $\psi_s \in [0,1)$ of the total biomass $B(t)$ to be made of structural components, e.g.,for a structural macromolecule ${P_s}$
\begin{align}\label{eq:biomass_composition_base}
\psi_s B(t) \leq {P_s}(t),~ \forall t\geq 0.
\end{align}
The extension of \eqref{eq:biomass_composition_base} to the network level can be expressed by collecting the individual constraints into the \emph{biomass composition matrix} $H_b$ with
\begin{align}\label{eq:biomass_composition}
H_b \begin{pmatrix}
C(t) \\ P(t)
\end{pmatrix} \leq 0,
\end{align}
where the rows of $H_b$ are derived from \eqref{eq:biomass_composition_base}.
We call \eqref{eq:biomass_composition} the \emph{biomass composition constraint}.
Furthermore, we can enforce specific reaction rates %~$v_m$
\begin{align}\label{eq:maintenance_base}
\begin{split}
 v_m(t) &\geq \phi_m B(t), ~ \forall t \geq 0\\
\end{split}
\end{align}
with the \emph{maintenance coefficient} $\phi_m \in [0,1)$ to model maintenance reactions scaling with biomass, e.g., re-synthesis
of lipids.
Hence, we call \eqref{eq:maintenance} the \emph{maintenance constraint}
\begin{align}\label{eq:maintenance}
v(t) \geq H_m \begin{pmatrix} C(t) \\ P(t) \end{pmatrix},
\end{align}
with the rows of $H_m$ corresponding to $\phi_m (w_C^T, w_P^T)$ (cf.~\eqref{eq:maintenance_base}). %or zero for non-maintenance reactions.
To construct the full deFBA problem, we introduce an \emph{end-time} $\tend>0$ and define the objective function as accumulation of the objective biomass \eqref{eq:objective_biomass} as
\begin{align}
\begin{split}\label{eq:deFBA_problem}
    \max_{v(t)} & \int_0^{\tend} B_o(t)\diff t \\
    \mathrm{s.t.}~&  \eqref{eq:qss},  \eqref{eq:ecc},\eqref{eq:box_constraints}, \eqref{eq:biomass_composition},\eqref{eq:maintenance}; \forall t \in [0,\tend].
\end{split}
\end{align}
This dynamic optimization problem can be solved by discretization with a collocation method. 
The result is a linear program (LP) for which efficient, specialized solvers are available.
% like \emph{gurobi} \citep{gurobi}, \emph{cvxopt} \citep{cvxopt}, or  SoPlex \citep{gleixner2016}.
With respect to the computational and numerical details of solving such problems, we refer the reader to \citep{waldherr2015}, and to \citep{reimers2017} for a large scale example.
We provide an implementation of the deFBA model class in Python 2.7\footnote{\url{https://bitbucket.org/hlindhor/defba-python-package}}, which imports/exports models using libSBML \citep{bornstein2008libsbml} and the resource allocation modeling (RAM) annotations \citep{ram2017}.
A step-by-step guide for the generation of deFBA models is described in \citep{reimers2017generating}.

\subsection{Important growth modes}
There are multiple reasons to discard the large end-time $\tend$ in favor of a shorter prediction horizon $0 < \tp << \tend$ and implement an iterative version of the original problem \eqref{eq:deFBA_problem}.
%Foremost,  the deFBA can produce artificial solutions, which do not reflect the behavior of the modeled organism.
%We call the artificial arcs \emph{linear phases}, which are defined by 
Foremost, the deFBA can produce \emph{linear phases}, defined as
\begin{align}
\diff B_o(t) / \diff t = \lambda,
\end{align}
with the constant linear growth rate $\lambda \geq 0$. 
These phases can occur if some macromolecules are very "cheap" in comparison to others. The model uses all resources to solely produce the cheap molecules, regardless of their utility. 
These phases can either be observed when using very small end-times or as mean to top off the objective value near nutrient depletion or the end-time $\tend$ \citep{waldherr2017}.
%thein near to nutrient depletion or near the end-time $t_{end}$ of the optimization horizon, to top off the objective value \citep{waldherr2017}.
We regard the linear phases as mathematical artifacts of the optimization method itself as we do not know of biological examples for this behavior.
Thus, one goal of the prediction horizon is to eliminate these linear arcs in the solutions.

Another important growth mode, called a \emph{balanced phase}, is defined by
\begin{align}
    \diff B_o(t) / \diff t  = \mu_{\bal} B_o(t) ,
\end{align}
with the  constant exponential growth rate $\mu_{\bal} \in \mathbb{R}_{\geq 0}$ depending on nutrient availability and the current biomass composition.
In these phases the composition of the biomass stays fixed as it is already optimal for the environment.
A dynamic solution generated by the deFBA typically consists of a series of balanced growth phases and the transitions between these.
%Based on our previous studies, solutions generated by the deFBA typically consist of balanced growth phases with adaptation phases in between and a short linear phase near the end-time \citep{waldherr2017}.
%Based on our previous studies, for most models and environmental conditions piece-wise balanced growth with a linear phase at the end seems to be the optimal solution.
%But we can not exclude more efficient exponential growth modes with varying $\mu(t)$, so we can use balanced phases as a lower bound to exponential growth.

\section{Short-term deFBA}
\subsection{Implementing the receding time horizon}\label{sec:receding_prediction_horizon}
The implementation of the receding prediction horizon $\tp$ is straightforward.
We split the time interval $[0, \tend]$ into intervals $[t_k, t_{k+1}]$ using the time grid $\Delta_t (\tc) = \{ t_k = k\tc~\vert~k \in \mathbb{N} \}$ defined by the \emph{iteration time} $\tc \in (0, \tp)$.
Then we replace the original deFBA problem \eqref{eq:deFBA_problem} with a series of small problems we call the \emph{short-term deFBA} (sdeFBA).
With given values ${Y}^{t_k}$, $C^{t_k}$, $P^{t_k}$, these read
\begin{subequations}\label{eq:sdeFBA_problem}
\begin{align}
\max_{v(t)} & \int_{t_k}^{t_{k}+\tp} B_o(t)\diff t \label{eq:sdeFBA_problem_objective}\\
\mathrm{s.t.}\; &\forall t \in [t_k,t_{k}+\tp] \\
& \frac{\mathrm{d}}{\mathrm{d}t} \begin{pmatrix} Y(t) \\ C(t) \\ P(t) \end{pmatrix} = \begin{pmatrix}S_{Y} \\ S_{C} \\ S_{P} \end{pmatrix} v(t) \label{eq:sdeFBA_problem_dynamics}\\
&  S_X v(t) = 0 \label{eq:sdeFBA_problem_steady_state}\\
& H_c v(t) \leq H_e P(t) \label{eq:sdeFBA_problem_ecc}\\
& H_b \begin{pmatrix} C(t) \\ P(t) \end{pmatrix} \leq 0 \label{eq:sdeFBA_problem_bcc}\\
& v(t) \geq H_m \begin{pmatrix} C(t) \\ P(t) \end{pmatrix}\label{eq:sdeFBA_problem_maintenance}\\
& v_{\min} \leq v(t) \leq v_{\max} \label{eq:sdeFBA_problem_box}\\
& Y(t),\;C(t),\;P(t)\geq 0 \label{eq:sdeFBA_problem_positivity}\\
& Y(t_k) = {Y}^{t_k}, ~ C(t_k)=C^{t_k}, ~ P(t_k) = P^{t_k}. \label{eq:sdeFBA_problem_continous}
\end{align}
\end{subequations}
For given initial values $Y_0, C_0, P_0$, we solve the problem iteratively starting at time zero and connecting the iterations via \eqref{eq:sdeFBA_problem_continous}.
The solution trajectories $Y^*(t)$, $C^*(t)$, $P^*(t)$, $v^*(t)$, ~$0 \leq t \leq \tend$ are generated by appending the calculated slices over the iteration time $[t_k , t_k + \tc]$ after each iteration.

\subsection{Choosing the prediction horizon}\label{sec:choosing_the_prediction_horizon}
\begin{figure}
    \centering
    \includegraphics[scale=0.7]{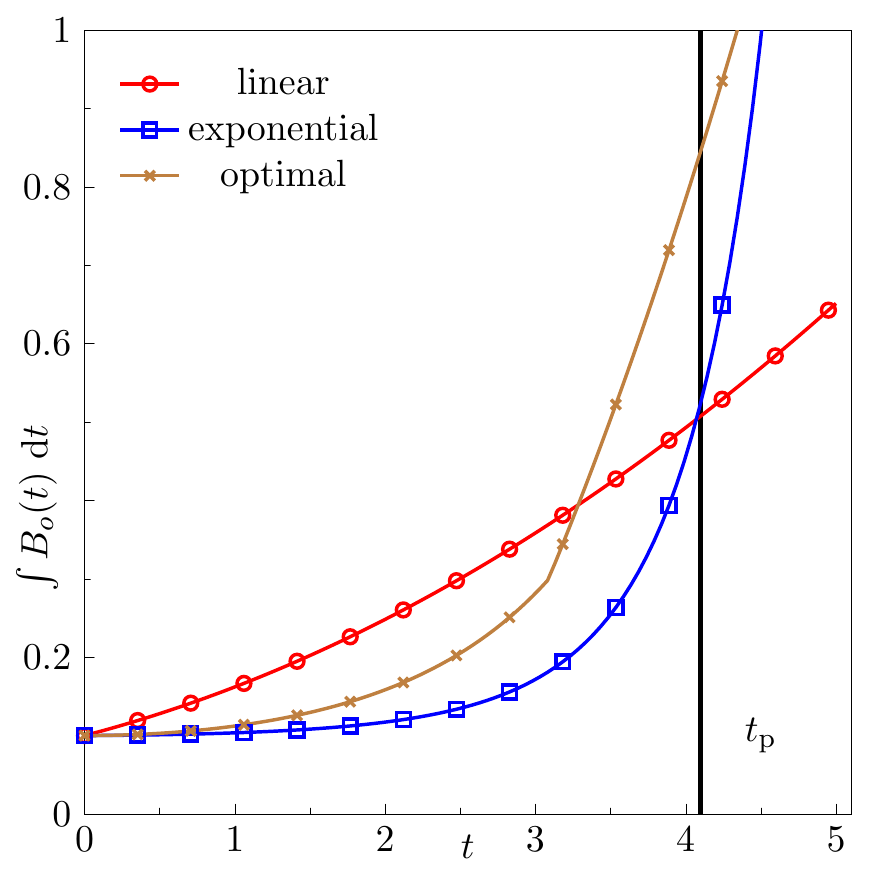}
    \caption{Illustration for choosing the prediction horizon. Upper bound on linear growth shown in red ($\circ$), balanced growth in blue ($\square$), and optimal solution in brown (x).}\label{fig:lin_vs_exp}.
\end{figure}
We already stated that the native growth mode for metabolic networks is exponential growth, while linear phases are undesired.
Our analysis in \citep{waldherr2017} shows, that linear solutions can arise on very short time scales as exponential solutions need a longer time horizon to outperform them.
Hence, we must ensure to choose the prediction horizon $t_p$ large enough such that linear solutions become sub-optimal.
At the same time we want to keep $t_p$ as small as possible to minimize computational cost.
We suggest to determine the prediction horizon by comparison of a strict upper bound on linear growth with an arbitrary balanced growth phase.
The idea is sketched in Figure \ref{fig:lin_vs_exp}.
This way we ensure the existence of at least piece-wise exponential solutions on the time horizon $\tp$.
This calculation is dependent on two sets of variables; the nutrients available and the initial biomass composition $P_\init, ~C_\init$ at time zero (or $t_k$ in the sdeFBA).
To eliminate the influence of nutrient availability in this first investigation we make the following assumption.
\begin{assumption}\label{as:nutrients}
    All external components $Y$ are limitlessly available.
\end{assumption}

We define the \emph{initial objective biomass} as 
\begin{align}
B_{\init} = b_C^T C_\init + b_P^T P_\init.
\end{align}

First we identify a strict upper bound on linear growth dependent on the initial biomass amount by constructing an optimization problem inspired by the regular FBA \citep{varma1994stoichiometric}.
We assume a linear growth phase $\diff B_o(t)/ \mathrm{d} t = \lambda$ and maximize the linear growth rate
\begin{align}
\lambda &= b^T_C S_C v_\lin + b^T_P S_P v_\lin.
\end{align}
Following Assumption \ref{as:nutrients}, we ignore the nutrient dynamics.
The optimization problem is then constructed as
\begin{subequations}\label{eq:upper_lin_bound}
\begin{align}
\lambda_s(B_{\init}) = \underset{\vlin,P_\lin, C_\lin}{\max}~~ & b^T_C S_{C} v_\lin + b^T_P S_{P} v_\lin \label{eq:upper_lin_bound_obj}\\
\mathrm{s.t.~~} & S_X \vlin = 0 \label{eq:upper_lin_bound_qss}\\
& H_c \vlin - H_e P_\lin \leq 0 \label{eq:upper_lin_bound_ecc}\\
& H_b \begin{pmatrix} C_\lin \\ P_\lin  \end{pmatrix} \leq 0 \label{eq:upper_lin_bound_bcc}\\
& w_C^T C_\lin + w_P^T P_\lin = B_\init \label{eq:upper_lin_bound_biomass}\\
& \vlin \geq H_m \begin{pmatrix} C_\lin \\ P_\lin  \end{pmatrix} \label{eq:upper_lin_bound_maintenance} \\
& v_{\min} \leq \vlin \leq v_{\max}, \label{eq:upper_lin_bound_box}
\end{align}
\end{subequations}
with \eqref{eq:upper_lin_bound_biomass} fixing the initial amount of biomass to $B_\init$.
The value of the \emph{specific growth rate} $\lambda_s(B_\init)$ is dependent on the amount of biomass.
Instead we use the \emph{regularized rate}
\begin{align}
\lambda_r = \frac{\lambda_s(B_\init)}{B_\init}.
\end{align}
For easier reading we omit the dependency of $\lambda_s$ on the biomass.

We construct the linear solution as 
\begin{align}
P(t) = P_\lin + S_P \vlin t,~ C(t)=C_\lin + S_C \vlin t.\label{eq:lin_solution}
\end{align}
This solution is usually not feasible for the original sdeFBA problem \eqref{eq:sdeFBA_problem} with $\tp>0$ as violations of \eqref{eq:sdeFBA_problem_bcc} and \eqref{eq:sdeFBA_problem_maintenance} are to be expected with increase in biomass over time.

As next step, we identify a balanced growth phase to use as a lower bound for optimal exponential growth by optimizing the static growth rate $\mu \geq 0$ at $t=0$
%To ensure exponential growth, we enforce balanced growth (cf. \citep{campbell1957synchronization}) at time zero
\begin{align}
\frac{\diff}{\diff t}\begin{pmatrix} C_\init \\ P_\init \end{pmatrix} = \mu \begin{pmatrix} C_\init \\ P_\init \end{pmatrix}.
\end{align}
%and maximize the growth rate .
The resulting optimization problem reads
\begin{subequations}
    \begin{align}
    \mu_{\bal} = & \max_{v_{\bal}}  \mu \\
     \text{s.t.~~}& \mu \begin{pmatrix} C_\init \\ P_\init \end{pmatrix} =\begin{pmatrix} S_C \\ S_P \end{pmatrix} v_{\bal} \\
     & S_X v_{\bal} = 0 \\ 
     & H_c v_{\bal} - H_E P_\init \leq 0 \\
     & v_{\bal} \geq H_m \begin{pmatrix} C_\init \\ P_\init  \end{pmatrix} \\
     &  v_{\min} \leq v_{\bal} \leq v_{\max}.
    \end{align}
\end{subequations}
%Contrary to the linear bound solution, the trajectories obtained by solving the initial value problem
The trajectories of the balanced growth phase are derived by solving the initial value problem
\begin{align}\label{eq:bal_ode}
    \frac{\diff}{\diff t}\begin{pmatrix} C(t) \\ P(t) \end{pmatrix} & = \mu_\bal \begin{pmatrix} C(t) \\ P(t) \end{pmatrix},
\end{align}
with $C(0)=C_\init$, $P(0)=P_\init$.
These trajectories are realized by the rates $v(t) = v_\bal  e^{\mu_{bal} t}$ and represent a feasible solution to \eqref{eq:sdeFBA_problem}, if Assumption \ref{as:nutrients} holds and the initial values are feasible
\begin{align}
 H_b \begin{pmatrix} C_\init \\ P_\init  \end{pmatrix} \leq 0.
\end{align}
% (cf.~\eqref{eq:sdeFBA_problem_bcc}) 

We can calculate a suitable time $\tp$, by comparing the the balanced solution \eqref{eq:bal_ode} to the linear one \eqref{eq:lin_solution}.
%With this we can calculate the time $\tp$ at which the balanced solution \eqref{eq:bal_ode} outgrows the upper bound on linear growth \eqref{eq:lin_solution}.
%We do this by comparing the integrated biomass curves, which correspond to the value of the objective function for both solutions over time.
The integral of the biomass curve for \eqref{eq:bal_ode} is derived as
\begin{align}
\begin{split}
IB_{\bal}(t, \mu_\bal, B_\init) & = \int_0^{t} b_C^T C(t) + b_P^T P(t) ~\mathrm{d}t \\
 &= \mu_{\bal}^{-1}  B_\init (e^{\mu_{\bal} {t}} - 1) \label{eq:exp_solution}
 \end{split}
\end{align}
and the corresponding integral for the linear case is 
\begin{align}
IB_{\lin}({t},\lambda_r,B_\init) & = \int_0^{t} B_\lin(t) ~\mathrm{d}t \\
                                                  & = \frac{\lambda_r B_\init}{2} {t}^2  + B_\init~ t. \label{eq:linear_solution}
\end{align}
We calculate the prediction horizon by solving 
\begin{align}\label{eq:determine_tp}
IB_{\mathrm{lin}}(\tp, \lambda_r,B_\init) - IB_{\bal}(\tp,\mu_{\bal}, B_\init)= 0
\end{align}
for $\tp$.
By looking at the slopes of the biomass curves at time zero, we can deduce that this $\tp >0$ only exists if, and only if, $\lambda_r > \mu_{\bal}$.
Otherwise, the model does not tend to the linear solution and we can chose $\tp$ arbitrarily.
\begin{assumption}\label{as:lambdageqmu}
The linear growth rate is larger than the balanced growth rate $\lambda_r > \mu_{\bal}$.
\end{assumption}
%The solution $\tp=0$ always exists for \eqref{eq:determine_tp}.
%If this is the only positive solution the curve $IB_{\bal}$ is larger for $t>0$.
%This means each solution of the sdeFBA on any time-scale grows exponentially and we can choose the prediction horizon arbitrarily.
%Secondly, there may exist another solution $\tp > 0$ for \eqref{eq:determine_tp}, at which the balanced growth solution outgrows the linear one (cf. Figure \ref{fig:lin_vs_exp}).
An optimal solution of \eqref{eq:sdeFBA_problem} on $[0, \tp]$ can only produce an objective value equal or larger than $IB_{\bal}(\tp)$, otherwise it would contradict the optimality principle.
Hence, we conclude that this optimal solution must contain a superlinear (typically exponential) arc as shown in Figure \ref{fig:lin_vs_exp}.

% In the next section we will show that this exponential phase must be at the start of the solution.
\begin{remark}\label{rm:recalculate}
Calculating $\tp$ is strongly dependent on the initial biomass $P_\init, C_\init$.
Hence, during an sdeFBA run the prediction horizon should be recalculated after each iteration step.
\end{remark}

\subsection{Choosing the iteration time}\label{sec:choosing_the_control_horizon}
To keep the computational cost of a sdeFBA run as small as possible we choose the iteration time $\tc$ as large as possible, such that the solution is still of exponential form.
Hence, we show that each solution of \eqref{eq:sdeFBA_problem} starts with an exponential phase.
For this we assume a solution starting with a linear phase
\begin{equation}\label{eq:B_mix}
 B_\mix(t) = \left\{ \begin{array}{ll}
     B_\init \lambda_r t +B_\init & 0 \leq t \leq \ts \\
     B_\init (\lambda_r \ts + e^{\mu_{\bal}(t - \ts)}) & \ts < t \leq \tp,
     \end{array} \right.
\end{equation}
with the switching time $\ts$ and assume Assumption \ref{as:lambdageqmu} holds.
This solution is constructed on the assumption that the linear growth phase does not benefit the autocatalytic capabilities of the system.
%We are using $\mu_{\max} \geq \mu(\ts)$ instead of a switching time dependent growth rate to avoid its inclusion to derivatives of $B_\mix$ with respect to $\ts$.
%Assumption \ref{as:nutrients} ensures $\mu_{\max}$ is unique and it can be identified by the RBA \citep{goelzer2011cell}.
\begin{theorem}
 If Assumption \ref{as:lambdageqmu} holds, any optimal solution curve $B_\mix$ \eqref{eq:B_mix} consists only of a single linear phase with $\ts = \tp$.
\end{theorem}
\begin{pf}
We identify the optimal switching time by solving
\begin{align}\label{eq:max_Bmix}
\max_{\ts} \int_0^{\tp} B_\mix(t) ~\diff t
\end{align}
analytically by finding local extrema via the first order derivative with respect to $\ts$
\begin{align}\label{eq:zero_theorem_1}
\begin{split}
      0 & = \frac{\diff} {\diff \ts} \int_0^{\tp} B_\mix(t) ~\diff t\\
        & = B_\init(\lambda_r(\tp - \ts) + 1 - e^{\mu_{\bal}(\tp-\ts)}),
\end{split}
\end{align}
with the obvious zero $\bar{t}_{\mathrm{s}} = \tp$.
Evaluating the second derivative at this point gives
%\begin{small}
%\begin{align}
%\begin{split}
%    \left.\frac{\diff^2} {\diff \ts^2} \int_0^{\tp} B_\mix (t) ~\diff t\right|_{\bar{t}_{\mathrm{s},1}} & = -\lambda_s + \mu_{\max}B_\init + \lambda_s \mu_{\max} \tp,  \\
%    \left.\frac{\diff^2} {\diff \ts^2} \int_0^{\tp} B_\mix (t) ~\diff t\right|_{\bar{t}_{\mathrm{s},2}} & = \lambda_s - \mu_{\max}B_\init - \lambda_s \mu_{\max} \tp.
%\end{split}
%\end{align}
%\end{small}
\begin{small}
\begin{align}
\begin{split}
    \left.\frac{\diff^2} {\diff \ts^2} \int_0^{\tp} B_\mix (t) ~\diff t\right|_{\bar{t}_{\mathrm{s}}} & = B_\init(\mu_{\bal}-\lambda_r) < 0,
\end{split}
\end{align}
\end{small}
\hspace*{-.2cm}with the last inequality following Assumption \ref{as:lambdageqmu}.
Hence, $\bar{t}_{\mathrm{s}}$ is a local maximum and any solution of the $B_\mix$ form does not include an exponential arc.
For the sake of completeness, we must also mention that there exists another zero of \eqref{eq:zero_theorem_1} $\bar{t}_{\mathrm{s,2}} \in [0,\tp)$, which cannot be given in closed form.
But, due to continuity and the intermediate value theorem, $\bar{t}_{\mathrm{s,2}}$ is a local minimum of \eqref{eq:max_Bmix}. \hfill~ \qed
%First assume 
%\begin{align}
%-\lambda_r + \mu_{\max} + \lambda_r \mu_{\max} \tp < 0.
%\end{align}
%Then $\bar{t}_{\mathrm{s},1}$ is a local maximum and $\ts = \tp$ holds.
%If on the other hand
%\begin{align}\label{eq:thm_1_condition}
% -\lambda_r + \mu_{\max} + \lambda_r \mu_{\max} \tp > 0 ~\Leftrightarrow \tp>\frac{1}{\mu_{\max}} - \frac{1}{\lambda_r},
%\end{align}
%holds, ${\bar{t}_{\mathrm{s},2}}$ is the local maximum.
%But \eqref{eq:thm_1_condition} also implies $\bar{t}_{\mathrm{s},2}>\tp$.
%Thus, in both cases the the switching time is larger than the prediction horizon and the whole solution is a single linear phase.
\end{pf}
As we have chosen $\tp$ such that the balanced growth solution \eqref{eq:exp_solution} outgrows the maximal linear one, we know that there exists a time frame $[0, \tc]$ on which the solution of \eqref{eq:sdeFBA_problem} must at least grow exponentially.
Thus, we assume the following form for the solution
\begin{small}\begin{align}\label{eq:B_opt}
    B_\opt(t) & = \left\{ \begin{array}{ll}
            B_\init e^{\mu_{\bal}t}, & 0 \leq t \leq \ts, \\
             B_\init  e^{\mu_{\bal}\ts} ( \lambda_r  (t-\ts) + 1), & \ts < t \leq \tp,
    \end{array}\right.
\end{align}\end{small}
\hspace*{-.2cm}with $B_\init \lambda_r e^{\mu_{\bal}\ts} = \lambda_s (B_\opt(\ts))$.
We want to choose $\tc$ such that no linear phase occurs in the final solution of the sdeFBA.
Otherwise, we can get faulty solutions as shown in the next section.

\begin{theorem}
    If Assumption \ref{as:lambdageqmu} holds, an optimal solution $B_\opt$ \eqref{eq:B_opt} of the sdeFBA \eqref{eq:sdeFBA_problem} is growing exponentially on the time frame $[0, \tc)$, with 
    \begin{align}\label{eq:approx_tc}
       0 < \tc  < \tp - 2\left( \frac{1}{\mu_\bal}- \frac{1}{\lambda_r}\right).
    \end{align}
\end{theorem}
\begin{pf}
    As in the previous proof, we identify the optimal switching time $\ts$ by solving the optimization problem
    \begin{align}\label{eq:theo_2_obj}
        \max_{\ts} \int_0^{\tp} B_\opt(t) ~\diff t.
    \end{align}
    The zeros of the first order derivative are given by
    \begin{gather}
%        \begin{split}
            \frac{\diff}{\diff \ts} \int_0^{\tp} B_\opt(t) ~\diff t = 0 \\
            \Rightarrow ~ \hat{t}_{\mathrm{s},1} = \tp - 2\left( \frac{1}{\mu_\bal}- \frac{1}{\lambda_r}\right),~ \hat{t}_{\mathrm{s},2} = \tp .
%        \end{split}
    \end{gather}
    The second-order derivative evaluated at these points is
    \begin{small}
    \begin{align}
    \begin{split}
    \left.\frac{\diff^2} {\diff \ts^2} \int_0^{\tp} B_\opt (t) ~\diff t\right|_{\hat{t}_{\mathrm{s},1}} & = (\mu_{\bal} - \lambda_r  ) B_\init e^{\mu_{\bal}\tp}< 0, \\
    \left.\frac{\diff^2} {\diff \ts^2} \int_0^{\tp} B_\opt§ (t) ~\diff t\right|_{\hat{t}_{\mathrm{s},2}} & = (\lambda_r - \mu_{\bal})B_\init e^{\mu_{\bal}\tp}> 0.
    \end{split}
    \end{align}
    \end{small}
    Hence, $\hat{t}_{\mathrm{s},1}$ maximizes \eqref{eq:theo_2_obj} and the solution is of exponential form until $\hat{t}_{\mathrm{s},1}$. \hfill~ \qed
\end{pf}
%Depending on the application of the sdeFBA it can be beneficial to choose the iteration time smaller than given by \eqref{eq:approx_tc}.
We strongly advise to choose the iteration time smaller than given by \eqref{eq:approx_tc} to compensate for numerical errors.
Otherwise, we might see solutions mixing linear and exponential phases as shown in Figure \ref{fig:results_1} (C).

Please note that $\tc$ is also dependent on the prediction horizon $\tp$ and the initial biomass composition $B_\init$.
So it should be recalculated together with $\tp$ after each iteration (cf. Remark \ref{rm:recalculate}).

\section{Numerical example}\label{sec:example_1}
We present a simple model, analyzed in detail in \citep{waldherr2017}, to give the reader an idea about the impact of end-times,  prediction horizons, and iteration times on the quality of the solution.
In this minimal example the organism can invest nutrients in either its' auto catalytic capabilities by investing in enzymes or it can produce non-catalytic components yielding a better nutrients-to-biomass ratio.  % $N$ in the production of catalyzing enzymes $E$ to increase its auto-catalytic capabilities or it can produce non-catalytic components $M$, which yield a better nutrients-to-biomass ratio.
%Additionally, an internal metabolite $A$ must be produced as interim step.
The three irreversible reactions of the network are
\begin{subequations}
    \begin{eqnarray}
    v_A  : & 1 ~N      & \rightarrow 1~ A \label{eq:example_1_network1} \\ 
    v_E  : & 1~N + 1~A & \rightarrow 1~ E \label{eq:example_1_network2} \\
    v_M  : & 1~N + 1~A & \rightarrow 1~ M. \label{eq:example_1_network3}
    \end{eqnarray}
\end{subequations}
%\textcolor{red}{explanation may be shortened.}
The external nutrient $N$ represents a collection of components necessary for growth, such as carbon, nitrogen, etc.
Further processed components made from these nutrients are collected as the internal metabolite $A$.
We differentiate the macromolecules into the group of enzymes $E$, collecting the whole enzymatic machinery needed for growth, and non-enzymatic macromolecules $M$.
These can be interpreted as storage components such as lipids, starch, or glycogen.

Assuming unlimited nutrients, we would expect a biological system to work exclusively in the exponential phase and produce no storage $M$ at all.
But the deFBA model \eqref{eq:deFBA_problem} may generate a solution containing linear phases depending on the system parameters and the end-time.

In this work we are only interested in the effects of the time variables and fix the system parameters to the values shown in Table \ref{tb:ex1_numerical_values}.
The numerical results using these values were all generated with our Python deFBA package\footnote{Available at \url{bitbucket.org/hlindhor/defba-python-package}} using a discretization step size $d = 0.1$ h and the initial values $E(0) = M(0) = 0.1$ mol.
\begin{table}
    \caption{Values used in the numerical example}\label{tb:ex1_numerical_values}
    \centering
    \begin{tabular}{lccccccc}
        \toprule
         $b_M~[\frac{\mathrm{g}}{\mathrm{mol}}]$ & $b_E~[\frac{\mathrm{g}}{\mathrm{mol}}]$ & $k_A~[\mathrm{h}^{-1}]$ & $k_M~[\mathrm{h}^{-1}]$ & $k_E~[\mathrm{h}^{-1}]$ \\
        \midrule
        15  & 10 & 1.5 & 2 & 1  \\
        \bottomrule
    \end{tabular}
\end{table}

\begin{figure}
    \centering
    \includegraphics[scale=0.89]{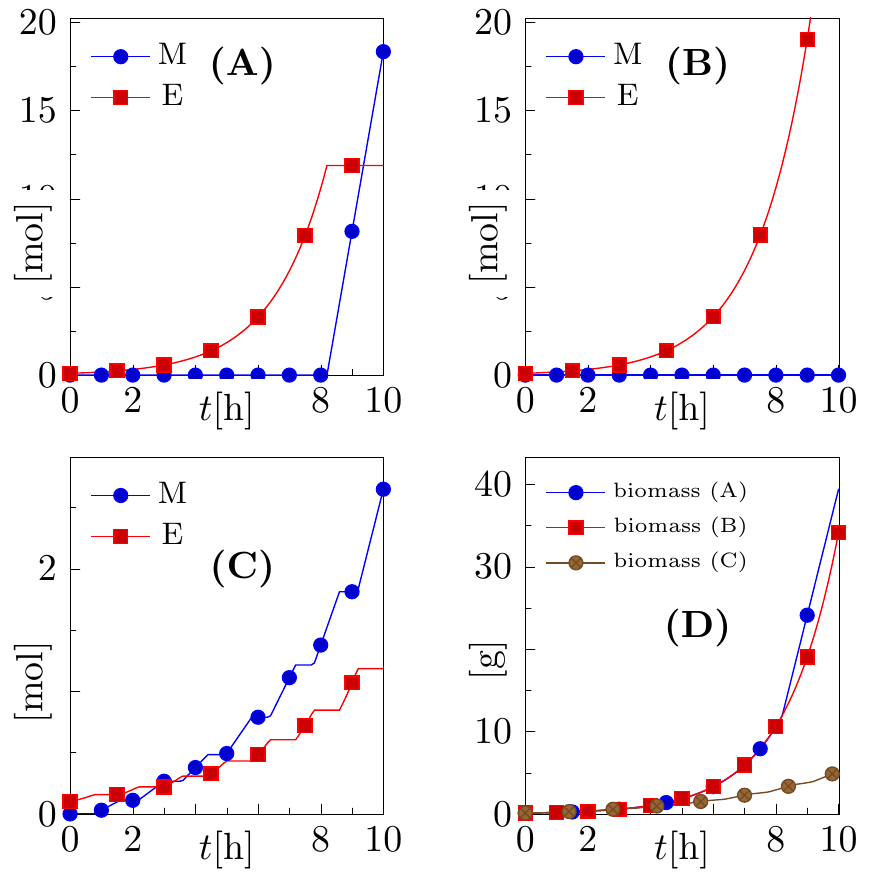}
    \caption{(A) deFBA solution $\tend=3$ h. (B) sdeFBA solution $\tp= 3.25$ h, $\tc=1.45$ h. (C) sdeFBA solution $\tp=2.5$ h, $\tc = 1.5$ h. (D) biomass comparison of methods. }\label{fig:results_1}
\end{figure}

Following \citep{waldherr2017}, we can derive the necessary condition for a single linear phase to be the optimal solution as
\begin{align}
t_{\mathrm{lin}} \leq \frac{2(k_M b_M - k_E b_E)}{b_M k_M k_E} \approx 1.45~\mathrm{h}.
\end{align}
Choosing any $\tend > t_{\mathrm{lin}}$ results in a mixed trajectory starting with an exponential phase and ending with a linear one.
This behavior can be observed in Figure \ref{fig:results_1} (A).
A purely exponential solution is not attainable with the deFBA as any solution ends in a linear phase producing only $M$ to top off the objective. 

But we can use the short-term deFBA to generate an exponential solution.
Using the idea from Section 3 we calculate the initial prediction horizon as $\tp \approx 3.25~\mathrm{h}$ and the iteration time as $\tc \approx 1.45~\mathrm{h}$.
The sdeFBA generates a purely exponential solution as shown in Figure \ref{fig:results_1} (B).
While this is a more reasonable solution from a biological view, the objective value for this solution is slightly smaller than the one obtained by the deFBA (cf. Figure \ref{fig:results_1} (D)).

Figure \ref{fig:results_1} (C) shows a sdeFBA solution using a prediction horizon $\tp = 2.5$ h and an iteration time $\tc = 1.5$ h. 
While this $\tp$ is capable of producing an exponential phase in each iteration the the chosen iteration time is way too large.
Hence, we see a solution in which exponential growth and linear phases take turns on each iteration slice.
This is neither optimal nor observed in nature.

\section{Conclusion}
While our presentation of the sdeFBA focuses on the quality of the solution, this method provides further advantages in comparison to the original deFBA.
Foremost, we can replace the fixed time frame $[0, \tend]$ in the original deFBA \eqref{eq:deFBA_problem} with a variable one dependent on the network's state.
As example, the deFBA is not designed to handle starvation scenarios and the optimization problem may become infeasible if the nutrients deplete.
But in the sdeFBA we can simply stop iterating once the nutrients deplete or another chosen threshold is reached.
Of course, this also means we can update state variables or dynamics while setting up the next iteration.
So we can use the sdeFBA as predictor in an online model predictive controller, which maximizes, e.g., some biomass component by changing the nutrient composition.

Lastly, the sdeFBA can be a way to solve large scale deFBA problems on large time-scales more efficiently.
The problem lies in the linear programs constructed by the deFBA, whose states can vary several orders of magnitude due to exponential growth phases. 
This leads to ill-posed problems, which take very long to solve even with sophisticated commercial solvers.
By breaking the problem into smaller pieces via the sdeFBA we can reduce the computational time.

\bibliography{Literatur.bib}
\end{document}